\documentclass{article}
\usepackage[american]{babel}
\usepackage{amsopn}
\usepackage{amsfonts}
\usepackage{hyperref}
\usepackage{amsmath,amssymb}
\usepackage{fancyhdr}
\numberwithin{equation}{section}

\def\xCn{{\rm C}}

\def\xR{{\mathbb R}}

\def\xN{{\mathbb N}}

\def\xS1{{\mathbb S}^1}

\def\xLtwo{{{\rm L}^2}}

\def\xHone{{{\rm H}^1}}

\def\xHn{{\rm H}}

\def\xLinfty{{\rm L}^{\infty}}

\newtheorem{thrm}{Theorem}[section]

\newtheorem{crllr}[thrm]{Corollary}
\newtheorem{prpstn}[thrm]{Proposition}

\newtheorem{rmrk}[thrm]{Remark}

\begin{document}
\title{Exit problems related to the persistence of solitons for the Korteweg-de Vries equation with small noise}

\author{{\sc By Anne de Bouard$^1$ and Eric Gautier$^2$}}


\maketitle

\begin{center}
$\ ^1$ Centre de Math\'ematiques Appliqu\'ees, UMR 7641 CNRS/Ecole
Polytechnique, 91128 Palaiseau cedex, France;
\texttt{debouard@cmapx.polytechnique.fr\vspace{0.2cm}}

$\ ^2$ ENSAE - CREST, 3 avenue Pierre Larousse, 92240 Malakoff,
France; \texttt{gautier@ensae.fr}
\end{center}

\begin{abstract}
We consider two exit problems for the Korteweg-de Vries equation
perturbed by an additive white in time and colored in space noise of
amplitude $\epsilon$. The initial datum gives rise to a soliton when
$\epsilon=0$. It has been proved recently that the solution remains
in a neighborhood of a randomly modulated soliton for times at least
of the order of $\epsilon^{-2}$. We prove exponential
upper and lower bounds for the small noise limit of the probability
that the exit time from a neighborhood
of this randomly modulated soliton is less than $T$,
of the same order in $\epsilon$ and $T$.
We obtain that the time scale is exactly the right one. We also study the
similar probability for the exit from a neighborhood of the deterministic
soliton solution. We are able to quantify the gain of eliminating the
secular modes to better describe
the persistence of the soliton.
\end{abstract}

\vspace{0.3cm}


\section{Introduction}\label{s1}

The Korteweg-de Vries (KdV) equation is a model for the evolution of
weakly nonlinear, shallow water, unidirectional long waves.
It is of the form
\begin{equation}
\label{KdV}
\partial_t u+\partial_x^3u+\partial_x (u^2)=0
\end{equation}
where the space variable $x$ is in $\xR$. The results of this paper
could easily be extended to generalized subcritical KdV equations for which the
nonlinearity is $\partial_x (u^p)$ for $p<5$, but we consider throughout
this article the $p=2$ case for simplicity. The KdV equation is
famous for its soliton solutions confirming the observation of the
solitary wave propagating on a channel by Russell in 1844. These
solitons are traveling waves of the form
$u_{c,x_0}(t,x)=\varphi_c(x-ct+x_0)$ where $c$ is the constant
velocity, $x_0\in\xR$ is the initial phase and
\begin{equation}
\varphi_c(x)=\frac{3c}{2\cosh^2\left(\sqrt{c}x/2\right)}
\end{equation}
These waves are localized, {\it i.e.} they decay exponentially to
zero as $x$ goes to infinity. Their shape is stable against
perturbations of the initial state. A first notion of stability, for
initial data close to $\varphi_{c}$, which takes into account the
symmetries of the evolution equation, is that of orbital stability.
This notion of stability
was first considered, for the solution $\varphi_c(x-ct+x_0)$ of the
KdV equation, by Benjamin \cite{Be}. The set $\{\varphi_c(\cdot-s),\
s\in\xR\}$ is the orbit of $\varphi_c$. The functional
$Q_c(u)=\mathbf{H}(u)+c\mathbf{M}(u)$ is used as a Lyapunov functional in the proof.
It involves two important invariant quantities of the evolution equation (\ref{KdV}):
the Hamiltonian, defined for $u$ in $\xHone(\xR)$, the space of square integrable functions
with square integrable first
order derivatives, by
\begin{equation}
\mathbf{H}\left(u\right)=\frac12\int_{\xR}\left(\partial_xu(x)\right)^2dx-\frac13\int_{\xR}u^3(x)dx
\end{equation}
and the mass defined by
\begin{equation}
\mathbf{M}\left(u\right)=\frac12\int_{\xR}u^2(x)dx.
\end{equation}
The space $\xHone(\xR)$ is the energy space, and it is a natural space for the solutions of (\ref{KdV}) :
indeed, if $u\in C([0,T]; \xHone(\xR))$ is a solution of (\ref{KdV}), then $\mathbf{H}(u(t))=\mathbf{H}(u_0)$
and $\mathbf{M}(u(t))=\mathbf{M}(u_0)$ for any $t\in [0,T]$, where $u_0$ is the initial datum in $\xHone$.
The shape of the soliton $\varphi_c$ is a solution of the constrained
variational problem which consists in minimizing the Hamiltonian for a constant mass.
Orbital stability means that when
the initial datum is close to $\varphi_c$ in $\xHone$ then the
solution remains close to the orbit of $\varphi_c$. The second stronger notion of stability is that
of asymptotic
stability. It states that, for initial datum close in
$\xHone$ to $\varphi_c$, the solution converges in some sense as
time goes to infinity to a soliton where the velocity and
phase have been shifted. Convergence may correspond to weak convergence in $\xHone$, see \cite{MM}
or \cite{MM2},
or strong convergence in some weighted Sobolev space, see \cite{PW}, for less
general perturbations of the initial datum,
when the solution is
written in the soliton reference frame. Note that strong convergence
in $\xHone$ is not expected due to the possibility of a dispersive
tail moving away from the soliton as time
goes to infinity.\\

\indent It is often physically relevant to consider random
perturbations of equation (\ref{KdV}), see \cite{SR}. It is also interesting
from a theoretical perspective
to study the stability of the soliton shape under these random perturbations.
We consider as in \cite{SR} the case of an additive noise which
could model a random pressure at the surface of the water. The
corresponding stochastic partial differential equation (SPDE for
short) written in It\^o form is the following
\begin{equation}\label{edps}
du+\left(\partial_x^3u+\partial_x(u^2)\right)dt=\epsilon dW
\end{equation}
where $\left(W(t)\right)_{t\ge0}$ is a Wiener process and $\epsilon$ is the small noise amplitude. As we
work in infinite dimensions, and in the absence of global smoothing property of the
group $S(t)$ on $\xHone$ associated to the unbounded operator $-\partial_x^3$,
$W$ needs to be a proper Wiener process on $\xHone$. Thus,
the components of $W(1)$ need to be
correlated for the law of $W(1)$ to be a {\it bona-fide} Gaussian
measure. It can then always be seen as the direct image via a
Hilbert-Schmidt self-adjoint mapping $\Phi$ of a cylindrical measure
and we assume that $\Phi$ is a mapping from $\xLtwo$ into $\xHone$. Recall that
$\Phi$ is Hilbert-Schmidt from $\xLtwo$ into $\xHone$ if it is a
bounded linear operator and for a complete orthonormal system
$\left(e_i\right)_{i\in\xN}$ of $\xLtwo$, $\sum_{i\in\xN}\left\|\Phi
e_i\right\|_{\xHone}^2$ is finite. The sum does not depend on the
complete orthonormal system and, endowed with its square root as a
norm, the space of Hilbert-Schmidt operators
$\mathcal{L}_2(\xLtwo,\xHone)=\mathcal{L}_2^{0,1}$ is a Hilbert
space. As a consequence, the Wiener process could be written as
\begin{equation}
W(t,x)=\sum_{i\in\xN}\beta_i(t)\Phi e_i(x),\quad t\ge0,\ x\in\xR
\end{equation}
where $\left(\beta_i\right)_{i\in\xN}$ is a collection of
independent standard real valued Brownian motions and
$\left(e_i\right)_{i\in\xN}$ a complete orthonormal system of
$\xLtwo$. Existence of path-wise mild solutions, almost surely
continuous in time for all $t$ positive with values in $\xHone$, of
the SPDE supplemented with the initial datum
$u(0)=u_0\in\xHone$ and uniqueness among those having almost surely
paths in some subspace $X_T\subset\xCn\left([0,T];\xHone\right)$ has
been obtained in \cite{dBD1}. In \cite{dBDT1,dBDT2} global well
posedness is obtained for rougher noises and less regular solutions.
It should be noticed that
in the physics literature the space-time white noise, corresponding to $\Phi=I$
is often considered, which we are not able to treat mathematically. For simplicity,
we consider the sequence of operators
\begin{equation}
\label{op}
\Phi_n=\left(I-\Delta+\frac1n(x^2I-\Delta)^{k}\right)^{-1/2},
\end{equation}
which are Hilbert-Schmidt
from $\xLtwo$ into $\xHone$ for $k$ large enough,
in order to prove lower bounds on exit times. As $n$ goes to infinity the Hilbert-Schmidt
assumption tends to be relaxed and the noise mimics a spatially homogeneous noise
with covariance $(I-\Delta)^{-1}$, which is a white noise in the Hilbert space 
$\xHone$.
It should also be noted that these operators are uniformly bounded in the space $\mathcal{L}^{0,1}$
of bounded operators from $\xLtwo$ into $\xHone$; indeed we have for every integer $n$,
$\|\Phi_n\|_{\mathcal{L}^{0,1}}\le 1$. It is possible to work with a more general approximating sequence,
see for example the kind of assumptions made in \cite{DG}.\\

\indent The linearized operator around
the soliton is particularly interesting to study the stability. It has a general null-space spanned by the two
secular modes $\partial_x\varphi_{c}$ and $\partial_{c}\varphi_c$. These modes are associated
with infinitesimal changes in the velocity and location of the solitary wave.
In Remark 2.3 in \cite{dBD2}, the following heuristic argument then implies that the random
solution should at most remain close to the deterministic solution up to times of the order
$\epsilon^{-2/3}$. It is based on an analogy with
the behavior of a linear system of SDEs such that 0 is a
degenerate simple eigenvalue corresponding to a Jordan block
$$\left\{
\begin{array}{l}
dX_1=X_2dt+\epsilon dW_1(t)\\
dX_2=\epsilon dW_2(t)
\end{array}\right.$$ with Brownian motions $W_1$ and $W_2$.
In such a case,
$$X_1(T)=\epsilon\int_0^TW_2(s)ds+\epsilon W_1(T)$$
has variance of the order of
$\epsilon^2T^3$ for large $T$. Thus, for a first approximation of the solution
$u^{\epsilon,\varphi_{c_0}}$ of (\ref{edps}) with initial datum $\varphi_{c_0}$, of the form
\begin{equation}\label{edet}
u^{\epsilon,\varphi_{c_0}}(t,x)=\varphi_{c_0}(x-c_0t)+\epsilon\tilde{\eta}^{\epsilon}(t,x-c_0t),
\end{equation}
and with an exit time defined for a neighborhood of the soliton in $\xHone$
$$B\left(\varphi_{c_0},\alpha\right)=\left\{f\in\xHone:\ \left\|f-\varphi_{c_0}\right\|_{\xHone}<\alpha\right\}$$
by
$$\tilde{\tau}_{\alpha}^{\epsilon}=\inf\left\{t\in[0,\infty):\
u^{\epsilon,\varphi_{c_0}}(t,\cdot+c_0t)\in
B\left(\varphi_{c_0},\alpha\right)^c \right\},$$ exit is expected to occur on a time scale of the order of $\epsilon^{-2/3}$.
However, it is believed that the soliton shape is preserved over a longer time scale.
A general approach which works for the deterministic equation, see \cite{MM,PW}, is to introduce a description by a
soliton ansatz where the parameters of the soliton fluctuate with time.
In the case of an additive noise physicists use an approximation of the
solution by a soliton ansatz of the form $\varphi_{c^{\epsilon}(t)}\left(x-x^{\epsilon}(t)\right)$
where
$c^{\epsilon}(t)$ and $x^{\epsilon}(t)$ are random scalar processes evolving according to a system of coupled SDEs.
In the case where the noise is the time derivative
of a one dimensional standard Brownian motion,
it is easily seen (see \cite{Wad}) that the solution can be written as
a modulated soliton plus a Brownian motion. However,
proving such a result is more involved when the noise is a function
of space as well. A mathematical justification is given in
\cite{dBD2} where the following result is proved.
\begin{thrm}\label{tdBD}
For $\epsilon>0$ and $c_0>0$, there exists $\alpha_0>0$ such that for every
$\alpha\in\left(0,\alpha_0\right]$ there exists a stopping time
$\tau_{\alpha}^{\epsilon}>0$ a.s. and semi-martingales
$c^{\epsilon}(t)$ and $x^{\epsilon}(t)$ defined a.s. for
$t\le\tau_{\alpha}^{\epsilon}$ with values in $(0,\infty)$ and $\xR$
respectively such that if we set
$$\epsilon\eta^{\epsilon}(t)=u^{\epsilon,\varphi_{c_0}}\left(t,\cdot+x^{\epsilon}(t)\right)-\varphi_{c^{\epsilon}(t)}$$
then
\begin{equation}\label{e1}
\int_{\xR}\eta^{\epsilon}(t,x)\varphi_{c_0}(x)dx=\left(\eta^{\epsilon},\varphi_{c_0}\right)=0,\quad\forall
t\le\tau_{\alpha}^{\epsilon}\quad a.s.,
\end{equation}
\begin{equation}\label{e2}
\int_{\xR}\eta^{\epsilon}(t,x)\partial_x\varphi_{c_0}(x)dx=\left(\eta^{\epsilon},\partial_x\varphi_{c_0}\right)=0,\quad\forall
t\le\tau_{\alpha}^{\epsilon}\quad a.s.
\end{equation}
and for all $t\le\tau_{\alpha}^{\epsilon}$,
\begin{equation}\label{e3}
\left\|\epsilon\eta^{\epsilon}(t)\right\|_{\xHone}\le\alpha
\end{equation}
and
\begin{equation}\label{e4}
\left|c^{\epsilon}(t)-c_0\right|\le\alpha.
\end{equation}
Moreover, there exists $C>0$ such that for all $T>0$ and
$\alpha\le\alpha_0$ there exists $\epsilon_0>0$ such that for all
$\epsilon<\epsilon_0$,
\begin{equation}\label{e5}
\mathbb{P}\left(\tau_{\alpha}^{\epsilon}\le
T\right)\le\frac{C\epsilon^2T\|\Phi\|_{\mathcal{L}_2^{0,1}}}{\alpha^4}.
\end{equation}
\end{thrm}
The proof uses the Lyapunov functional $Q_c$ as a central tool.
The equations \eqref{e1} and \eqref{e2} are such that restricted to
this subspace, the Lyapunov functional is coercive, {\it i.e.} the
operator $Q_{c_0}''$ is positive.
It allows to keep $|c^{\epsilon}(t)-c_0|$ and
$\|\epsilon\eta^{\epsilon}(t)\|_{\xHone}$ small on a longer time interval.
Also, these two conditions
together with the implicit function theorem allow to obtain $c^{\epsilon}$ and
$x^{\epsilon}$. Other results in \cite{dBD2} give
the asymptotic distribution of $\eta^{\epsilon}$ as $\epsilon$ goes
to zero as well as coupled equations for the
evolution of the random scalar parameters.
The parameters and remainder do not
depend on $\alpha \le \alpha_0$. In the upper bound \eqref{e5} the
product $\epsilon^2T$ appears. The theorem says that the solution
stays in a neighborhood of the randomly modulated soliton
$\varphi_{c^{\epsilon}(t)}(x-x^{\epsilon}(t))$ with high probability
at least for times small compared to $\epsilon^{-2}$. The time spent in
a neighborhood of a soliton-like wave, when the initial datum gives
rise to a soliton for $\epsilon=0$, is called the persistence time,
see \cite{DP} for numerical confirmations that the above order is
the right order.\\

\indent In this paper we first study the exit time $\tilde{\tau}_{\alpha}^{\epsilon}$
and obtain that the time scale on which the solution stays close to the deterministic
soliton is indeed at most of the order of $\epsilon^{-2/3}$. We then revisit the
upper bound \eqref{e5} and prove a sharper exponential bound. This bound
is supplemented with an exponential lower bound of the same order in the parameters $T$ and $\epsilon$.
We thus obtain the right order of the cumulative distribution
function (CDF) of the exit time off neighborhoods of the randomly
modulated soliton. This gives a confirmation that the time scale on which the
approximation of the solution by a randomly modulated soliton is valid is of the order $\epsilon^{-2}$.
Our main tools are large deviations along with a study of the associated variational problems.
Similarly, factors $T$ and $T^3$ have also
been obtained in the study of the tails of the mass and arrival time
for stochastic nonlinear Schr\"odinger equations (NLS) in \cite{DG,EG1} with the same techniques.
These quantities are the main processes impairing soliton transmission in optical fibers.
In that setting, physicists again use the approximation by a randomly modulated
soliton. An analogue of Theorem \ref{tdBD} for stochastic NLS
equations would allow to tell up to what length of the fiber line the approximation is licit.
Large deviations are also known to be a
useful tool to study the exit problem from an asymptotic equilibrium
point or noised induced transition between several equilibrium
points in the small noise limit (see \cite{FW}, and \cite{EG4} for an
exit problem for stochastic weakly damped nonlinear Schr\"odinger
equation). Here we however study a simpler problem than the escape
from the asymptotically stable central manifold
for the KdV equation which we hope to study in future
works.

\section{Large deviations and escape from a neighborhood of the soliton}
We use sample path large deviations in this article in order to obtain lower bounds
of the asymptotic as $n$ goes to infinity and $\epsilon$ goes to zero of
probabilities $\mathbb{P}\left(\tau^{n,\epsilon}\le T\right)$ where $\tau^{n,\epsilon}$
is the exit time of a neighborhood of either the deterministic soliton or the randomly modulated soliton.
The $n$ recalls that we consider a
sequence of operators $\Phi_n$, see (\ref{op}).
Large deviation techniques, see for example \cite{DZ,DS1},
allow to quantify convergence to zero of rare events.
For example, it is easy to check that on a finite time interval $[0,T]$, the paths of the
solution of
\ref{edps} starting form $\varphi_{c_0}$ converge in probability to the paths of the deterministic
soliton solution. The probability that exit from a neighborhood of the soliton occurs before $T$
goes to zero as $\epsilon$ goes to zero. Large deviations quantify the convergence to zero
of such probabilities. Following Varadhan's formalism, large deviations could be stated
as a sequence of inequalities called a large deviation principle (LDP for short).
The convergence to zero of the logarithm of the probabilities of rare events is characterized by
a speed, here $\epsilon^2$, and a deterministic functional $I^n$ depending on the operator $\Phi_n$ considered,
called rate function, to be minimized
on the closure and interior of the set defining the rare event in the state space.
In the small noise asymptotics and for sample path large deviations,
the rate function could be expressed
in terms of the mild solution of the control problem
\begin{equation}
\label{econtrol} \left\{\begin{array}{l} \partial_t u+\partial_x^3u+\partial_x (u^2)=\Phi_n h,\\
u(0)=\varphi_{c_0}\ \mbox{and}\ h\in
\xLtwo\left(0,T;\xLtwo\right).\end{array}\right.\end{equation} We
denote the solution by $\mathbf{S}^{n,\varphi_{c_0}}(h)$. The
mapping $h\to \mathbf{S}^{n,\varphi_{c_0}}(h)$ is called the control
map and \eqref{econtrol} the control equation. We recall that a rate
function $I$ on the sample space (here the paths space
$\xCn([0,T];\xHone)$) is lower semicontinuous and that a good rate
function is such that $I^{-1}([0,R])$ is compact for every $R$
positive.
\begin{thrm}\label{t1}
The laws $\left(\mu^{u^{n,\epsilon,\varphi_{c_0}}}\right)_{\epsilon>0}$
of the paths of the solutions of \eqref{edps} for the operator $\Phi=\Phi_n$ on
$\xCn([0,T];\xHone)$ with initial datum $\varphi_{c_0}$
satisfy a LDP of speed $\epsilon^2$ and good rate function
\begin{equation*}
I^n(w)=\frac12\inf_{h\in\xLtwo(0,T;\xLtwo):\
w=\mathbf{S}^{n,\varphi_{c_0}}(h)}\|h\|
_{\xLtwo\left(0,T;\xLtwo\right)}^2.
\end{equation*}
It means that for every Borel set $B$ of
$\xCn\left([0,T];\xHone\right)$, we have the lower bound
\begin{equation*}
-\inf_{{\tiny \begin{array}{rc}\circ\\w\in
B\end{array}}}I^n(w)\leq\underline{\lim}_{\epsilon\rightarrow0}\epsilon^2\log
\mathbb{P}\left(u^{n,\epsilon,\varphi_{c_0}}\in B\right)
\end{equation*}
and the upper bound
\begin{equation*}
\overline{\lim}_{\epsilon\rightarrow0}\epsilon^2\log\mathbb{P}\left(u^{n,\epsilon,\varphi_{c_0}}\in
B\right)\leq-\inf_{w\in\overline{B}}I^n(w).
\end{equation*}
\end{thrm}

The proof uses the LDP for the laws of the stochastic convolution
$\epsilon Z$ where $Z(t)=\int_0^{t}S(t-s)dW(s)$ on the Banach path space
$X_T$; it is a subspace of
$\xCn\left([0,T];\xHone\right)$ where the fixed point argument proving
the local well-posedness is used, see \cite{dBD1}.
The stochastic convolution appears when we write the equation satisfied by the mild solution
of \eqref{edps}. These laws are Gaussian measures and the LDP is a consequence of the general result on
LDP for centered Gaussian measures on real Banach spaces, see
\cite{DS1}. The second step is to prove the continuity of the mapping which, to the
perturbation $Z$ in $X_T$ assigns the
solution $u^{n,1,\varphi_{c_0}}:=\mathcal{G}(Z)$ in $\xCn\left([0,T];\xHone\right)$.
It is obtained noting that $\mathcal{G}(Z)=v(Z)+Z$ where $v(x)$ denotes the solution of
$$\left\{\begin{array}{l}
    \partial_t v+\partial_x^3v+\partial_x\left((v+Z)^2\right)=0\\
    v(0)=\varphi_{c_0}.
  \end{array}\right.$$
Then the continuity of $\mathcal{G}$ is a consequence of the continuity of $v$
with respect to the perturbation $Z$. It could be proved as in \cite{dBD3}
where the stochastic NLS equation is considered. LDP for the paths of the mild solution
of the SPDE is then obtained by a direct application of the contraction principle which
states that we can push forward LDP for measures on a Hausdorff topological space to
a LDP for direct image measures on another Hausdorff
topological space when the mapping is continuous. More details on the proof of such LDP are
given in \cite{EG1} where the stochastic NLS equation with additive noise is considered.
The control map  $\mathbf{S}^{H,\varphi_{c_0}}(h)$ for the theoretical stochastic equation with spatially
homogeneous noise
is defined as the mild solution of
\begin{equation}
\label{econtrolWN} \left\{\begin{array}{l} \partial_t u+\partial_x^3u+\partial_x (u^2)=\Phi_{H} h,\\
u(0)=\varphi_{c_0}\ \mbox{and}\ h\in
\xLtwo\left(0,T;\xLtwo\right),\end{array}\right.
\end{equation} where $\Phi_{H}=(I-\Delta)^{-1/2}$. Note that though we cannot give a mathematical
meaning to the stochastic equation with such homogeneous noise, the corresponding control map is well defined.\\

Let us now consider the exit times $\tilde{\tau}_{\alpha}^{n,\epsilon}$. Note that we only consider
here a lower bound of the probability since it is enough to prove the heuristic of Remark 2.3 in \cite{dBD2}.
Recall that we want to check that the time scale on which the approximation by the
deterministic soliton is licit is at most of the order of
$\epsilon^{-2/3}$. When studying the exit times $\tau_{\alpha}^{n,\epsilon}$, however, we give both upper
and lower bounds of the tail probabilities of the same order in the parameters $\epsilon$ and $T$.

\begin{prpstn}\label{p3}
Take $T$, $c_0$ positive; then for $\alpha_0$ small enough, for every
$\alpha<\alpha_0$, there
exists a constant $C(\alpha,c_0)$ which depends on $\alpha$ and
$c_0$ but not on $T$ such that
$$\underline{\lim}_{n\rightarrow\infty}\underline{\lim}_{\epsilon\to0}\epsilon^2\log\mathbb{P}
\left(\tilde{\tau}_{\alpha}^{n,\epsilon} \le
T\right)\ge-\frac{C(\alpha,c_0)}{T^3}.$$
\end{prpstn}

\noindent
{\bf Proof.} For fixed $n$, using Theorem \ref{t1},
$\underline{\lim}_{\epsilon\to0}\epsilon\log\mathbb{P}
\left(\tilde{\tau}_{\alpha}^{n,\epsilon}\le
T\right)$ is larger than
$$-\frac12\inf\left\{\|h\|_{\xLtwo(0,T;\xLtwo)}^2,\ h:\
\left\|\mathbf{S}^{n,\varphi_{c_0}}(h)(T)-\varphi_{c_0}(\cdot-c_0T)\right\|_{\xHone}
>\alpha\right\}.$$
In a first step, we consider the preceding variational problem in which the operator
$\Phi_n$ is replaced by the operator $\Phi_{H}=(I-\Delta)^{-1/2}$ and $\alpha$ is replaced by
$2\alpha$.
We give upper bounds on the infimum by minimizing on smaller and smaller sets of controls,
until we are able to handle the variational problem. Note that up to now, the problem
is a control problem for the KdV equation that we cannot handle. We will show that we
can work on more restrictive classes of controls and still obtain a nice qualitative order.
Using the Sobolev embedding of $\xHone$ into $\xLinfty$, with norm
$C_{\infty}$, the infimum is found to be less than
$$\inf\left\{\|h\|_{\xLtwo(0,T;\xLtwo)}^2,\ h:\
\left\|\mathbf{S}^{H,\varphi_{c_0}}(h)(T)-\varphi_{c_0}(\cdot-c_0T)\right\|_{\xLinfty}
>2C_{\infty}\alpha\right\}$$
since we then minimize on a smaller set.
We consider controls $h$ giving rise to modulated solitons of the form
$$\varphi_{c(t)}\left(x-\int_0^{t}c(s)\right)$$ in the homogeneous case. They are such that
$c(0)=c_0$ and
\begin{equation*}\Phi_{H}
h(t,x)=c'(t)\left.\partial_c\varphi_c\right|_{c=c(t)}\left(x-\int_0^tc(s)ds\right),
\end{equation*}
since the soliton profile $\varphi_c$ satisfies the equation
$$
-c \partial_x \varphi_c +\partial_x^3 \varphi_c +\partial_x (\varphi_c)^2=0.
$$
Again, taking the infimum on a smaller set of controls, we obtain the lower bound
\begin{align*}
-\frac12\inf&\left\{\int_0^T\left\|c'(t)\Phi_{H}^{-1}\left[
\left.\partial_c\varphi_c\right|_{c=c(t)}\left(x-\int_0^tc(s)ds\right)\right]
\right\|_{\xLtwo}^2dt,\right.\\
&\quad\left.c\in\xCn^1([0,T]; (0,+\infty)):\ c(0)=c_0,\
\left|\varphi_{c(T)}\left(\int_0^T(c_0-c(s))ds\right)-\frac{3c_0}{2}
\right|>2C_{\infty}\alpha\right\},
\end{align*}
where we have bounded from below the $\xLinfty$ norm of the function by its value at
$c_0T$.
This is in turn bigger than
\begin{align*}
&-\frac12\inf\left\{\int_0^T(c'(t))^2\left\|\Phi_{H}^{-1}
\left(\left.\partial_c\varphi_c\right|_{c=c(t)}\right)
\right\|_{\xLtwo}^2dt,\right.\\
&\quad\left.c\in\xCn^1([0,T]; (0,+\infty)):\ c(0)=c_0,\
\frac{3c_0}{2}-\varphi_{c(T)}\left(\int_0^T(c_0-c(s))ds\right)>2C_{\infty}\alpha\right\},
\end{align*}
due to the fact that $\Phi_H$ commutes with spatial translations.
Let us fix $\alpha_0<\frac{3c_0}{4C_{\infty}}$ so that
$c_0-\frac{4C_{\infty}\alpha}{3}>0$ for $0<\alpha<\alpha_0$.
A sufficient condition for the constraint on the terminal value to hold is
$$c_0-\frac{4C_{\infty}\alpha}{3}>\frac{4c(T)}{\exp\left(\sqrt{c(T)}\int_0^T(c_0-c(s))ds\right)}.$$
Noticing that the function $\lambda$ defined by
$\lambda(x)=4x^2\exp\left(-x\int_0^T(c_0-c(s))ds\right)$ attains its
maximum at $x=2/\int_0^T(c_0-c(s))ds$ for $x\ge0$, if $\int_0^T(c_0-c(s))ds\ge 0$,
we obtain that it is enough to
have
$$\int_0^T(c_0-c(s))ds>\frac{4}{e\sqrt{c_0-4C_{\infty}\alpha/3}} := \delta(c_0,\alpha).$$
As in \cite{DG} for the tails of the arrival time of a pulse driven by a
stochastic nonlinear Schr\"odinger equation where we obtained the order $-CT^{-3}$,
the boundary condition is in integrated form.
The integral to be minimized is of the form
$$\int_0^T\left(c'(t)\right)^2g\left(c(t)\right)dt$$ where
$g(c)=\|(I-\Delta)^{1/2}
\partial_c\varphi_c\|_{\xLtwo}^2$. Instead of solving the
problem of the calculus of variations with a nonstandard boundary
condition, we make a guess and look for solutions of the form
$c(t)=c_0-2\gamma t/T^2$ for some positive $\gamma$ with $c_0-2\gamma/T>0$. Note that if
$\gamma=\inf\{\frac32 \delta(c_0,\alpha),\frac{c_0}{4}\}$, then the boundary
condition is satisfied for $T\ge 1$. Also the term $g(c(t))$ in the integral is then
such that there exists a constant $C(c_0)$ with
$$\int_0^T\left(c'(t)\right)^2g\left(c(t)\right)dt \le C(c_0)\int_0^T\left(c'(t)\right)^2dt, \; \mathrm{for}
\; T\ge 1,$$
since $c_0/2\le c(t) \le c_0$ for any $T\in [0,T]$.
Thus, for a new constant $C(\alpha,c_0)$, we obtain
$$\int_0^T\left(c'(t)\right)^2g\left(c(t)\right)dt \le \frac{C(\alpha,c_0)}{T^3},\; \mathrm{for}
\; T\ge 1. $$
Let us now consider the case where the square root of the covariance operator of the noise is $\Phi_n$,
and let us start from the path $c(t)$ exhibited in the homogeneous noise case; we denote the corresponding
control by
$$h_c(t,x)=c'(t)\left(I-\Delta\right)^{1/2}
\left(\left.\partial_c\varphi_c\right|_{c=c(t)}\left(\cdot-\int_0^tc(s)ds
\right)\right).$$ Then, since for such $h_c$
$$\left\|\mathbf{S}^{H,\varphi_{c_0}}(h_c)(T)-\varphi_{c_0}(\cdot-c_0T)\right\|_{\xHone}>2\alpha,$$
we deduce from the continuity of the mild-solution
of the control map with respect to the convolution of the semi-group
with the control, used to prove the LDP, and the continuity
of this last quantity with respect to the control, that for sufficiently large $n$
\begin{align*}
&\left\|\mathbf{S}^{n,\varphi_{c_0}}(h_c)(T)-\varphi_{c_0}(\cdot -c_0T)\right\|_{\xHone}\\
&=\left\|\mathbf{S}^{H,\varphi_{c_0}}\left((I-\Delta)^{1/2}(I-\Delta+\frac{1}{n}(x^2I-\Delta)^k)^{-1/2}h_c\right)(T)
-\varphi_{c_0}(\cdot -c_0T)\right\|_{\xHone}\\
&>\alpha.
\end{align*}
This ends the proof.\hfill $\square$\\

\indent
As a consequence, the time scale on which an exit from a neighborhood of the
soliton occurs is at most $\epsilon^{-2/3}$. In the next section we prove that the time scale on which the
solution remains close to the randomly modulated soliton is exactly of the much longer order $\epsilon^{-2}$.
We provide upper and lower bounds for this result.

\section{Escape from a neighborhood of the randomly modulated soliton}
In the proof of Theorem 2.1 in \cite{dBD2} a local parametrization $u\mapsto(\mathcal{C}(u),\mathcal{X}(u))$
is used, in order to obtain parameters of the soliton wave form such that
$u=\varphi_{\mathcal{C}(u)}(\cdot-\mathcal{X}(u))
+\mathcal{R}(u)$ with $\mathcal{R}$ satisfying some adequate orthogonality conditions.
This parametrization is obtained using the implicit function theorem, imposing that the constraints
(\ref{e1}) and (\ref{e2}) hold. Such a parametrization holds as long as $u$ remains in a proper
neighborhood of the spatial translates of $\varphi_{c_0}$; thus, setting $c^{\epsilon}(t)
= \mathcal{C}(u^{\epsilon,\varphi_{c_0}}(t))$ and $x^{\epsilon}(t)=\mathcal{X}(u^{\epsilon,\varphi_{c_0}}(t))$,
with $u^{\epsilon,\varphi_{c_0}}$ a
solution of (\ref{edps}) with paths a.s. in $\xCn^1(\xR^+;\xHone)$ and initial datum $\varphi_{c_0}$,
the processes $c^{\epsilon}(t)$ and $x^{\epsilon}(t)$ are well defined adapted processes,
up to a stopping time of the form
$$\overline{\tau}_{\alpha}^{\epsilon}=\inf\left\{t\ge0,\ |c^{\epsilon}(t)-c_0|\ge \alpha\ {\rm or}\ \|u^{\epsilon,\varphi_{c_0}}(t,\cdot+x^{\epsilon}(t))-\varphi_{c_0}\|_{\xHone}\ge\alpha\right\}.$$
We can indeed always replace a function $u$ by the function $u(\cdot +\mathcal{X}(u))$
and come back to the case where $x^{\epsilon}(t)$ is close to 0, and it is not necessary
to include a condition of the form $|x^{\epsilon}(t)|\ge\alpha$. It is also shown that
this stopping time could be bounded above and below by
$$\tau_{C\alpha}^{\epsilon}=\inf\left\{t\ge0,\ |c^{\epsilon}(t)-c_0|\ge \alpha \ {\rm or}\ \|u^{\epsilon,\varphi_{c_0}}(t,\cdot+x^{\epsilon}(t))-\varphi_{c^{\epsilon}(t)}\|_{\xHone}\ge\alpha\right\}$$
for some constants $C$ depending solely on $c_0$ and $\alpha_0$, with $0<\alpha\le \alpha_0$.
Hence, the two stopping times are equivalent and the qualitative behavior of
$\mathbb{P}(\tau_{\alpha}^{\epsilon}\le T)$ with respect
to $\epsilon$ and $T$ should be the same as the original $\mathbb{P}(\overline{\tau}_{\alpha}^{\epsilon}\le T)$.
These new stopping times $\tau_{\alpha}^{\epsilon}$ prove to be more convenient to work with.
Note that the implicit function theorem says that each function in $\xHone$ sufficiently close to a translate of the soliton
could be written as a translated soliton with slightly different parameters plus a remainder
in the subspace of $\xHone$ orthogonal in $\xLtwo$ to $\varphi_{c_0}$ and $\partial_x\varphi_{c_0}$.
Therefore the exit in terms of the $\xHone$ norm is the exit of a proper open subset.\\

\indent Let us first revisit the upper bound given in \cite{dBD2} and prove that
the upper bound is indeed exponential. The proof relies on
exponential tail estimates. We denote by $\mathcal{L}^{-1,0}$ and $\mathcal{L}^{0,1}$ the
spaces of bounded operators from $\xHn^{-1}$ to $\xLtwo$ (respectively from $\xLtwo$ to $\xHone$).

\begin{prpstn}\label{p1}
For $T>0$ and $0<\alpha\le\alpha_0$ fixed,
there exists a constant $C(\alpha,c_0)$, depending on $\alpha$ and $c_0$ but not on $T$, and
$\epsilon_0>0$ with $\epsilon_0^2T$ sufficiently small depending on
$\|\Phi\|_{\mathcal{L}_2^{0,1}}$ and $\alpha$, such that for every positive
$\epsilon<\epsilon_0$,
\begin{equation}\label{e6}
\mathbb{P}\left(\tau_{\alpha}^{\epsilon}\le
T\right)\le\exp\left(-\frac{C(\alpha,c_0)}{\epsilon^2T\|\Phi\|_{\mathcal{L}^{0,1}}^2}\right).
\end{equation}
\end{prpstn}

\noindent
{\bf Proof.} Fix $T$ positive. The estimate \eqref{e5} relies on the
two following inequalities. Let
$\tau=\tau_{\alpha}^{\epsilon}\wedge T$; then for $\alpha_0$ sufficiently small there exists a positive constant
$C$ independent of $T$ such that
\begin{equation*}
|c^{\epsilon}(\tau)-c_0|^2\le
C\left[\|\epsilon\eta^{\epsilon}(\tau)\|_{\xLtwo}^4+4\epsilon^2\left|\int_0^{\tau}
(u^{\epsilon,\varphi_{c_0}}(s),dW(s))_{\xLtwo}\right|^2+\epsilon^4
\tau^2
\|\Phi\|_{\mathcal{L}_2^{0,1}}^4\right]
\end{equation*}
and
\begin{align*}
\|\epsilon\eta^{\epsilon}(\tau)\|_{\xHone}^2\le&
C\left[\|\epsilon\eta^{\epsilon}(\tau)\|_{\xLtwo}^4+4\epsilon^2\left|\int_0^{\tau}
(u^{\epsilon,\varphi_{c_0}}(s),dW(s))_{\xLtwo}\right|^2+\epsilon^4\tau^2\|\Phi\|_{\mathcal{L}_2^{0,1}}^4\right.\\
&\quad+\epsilon\int_0^{\tau}\left(\partial_xu^{\epsilon,\varphi_{c_0}}(s),\partial_xdW(s)\right)_{\xLtwo}
-\epsilon\int_0^{\tau}\left(\left(u^{\epsilon,\varphi_{c_0}}(s)\right)^2,dW(s)\right)_{\xLtwo}\\
&\quad+c_0\epsilon\int_0^{\tau}\left(u^{\epsilon,\varphi_{c_0}}(s),dW(s)\right)_{\xLtwo}
+\frac{\epsilon^2}{2}\tau\|\Phi\|_{\mathcal{L}_2^{0,1}}^2\\
&\quad\left.+\epsilon^2\|\Phi\|_{\mathcal{L}_2^{0,1}}^2\int_0^{\tau}\|u^{\epsilon,\varphi_{c_0}}(s)\|_{\xLtwo}ds
+c_0\frac{\epsilon^2}{2}\tau\|\Phi\|_{\mathcal{L}_2^{0,1}}^2\right].
\end{align*}
These are obtained from several manipulations of the
Lyapunov functional and evolution equations for the mass and
Hamiltonian evaluated on the solution of the stochastic KdV equation (see \cite{dBD2}).
The evolution of these quantities is obtained using the It\^o formula and a smoothing procedure.
We do not reproduce the proof here.
We may also write
\begin{align*}
\mathbb{P}\left(\tau_{\alpha}^{\epsilon}\le
T\right)&\le\mathbb{P}\left(|c^{\epsilon}\left(\tau\right)-c_0|^2\ge\alpha^2\
{\rm or}\ \|\epsilon\eta^{\epsilon}\left(\tau\right)\|_{\xHone}^2\ge\alpha^2\right)\\
&\le\mathbb{P}\left(|c^{\epsilon}\left(\tau\right)-c_0|^2\ge\alpha^2
\right)+\mathbb{P}\left(\|\epsilon\eta^{\epsilon}\left(\tau\right)\|_{\xHone}^2\ge\alpha^2\right).
\end{align*}
Note that when $|c^{\epsilon}\left(\tau\right)-c_0|^2\ge\alpha^2$ we have
$|c^{\epsilon}\left(\tau\right)-c_0|^2=\alpha^2$ and
$\|\epsilon\eta^{\epsilon}\left(\tau\right)\|_{\xLtwo}\le\alpha$. Thus for
$\epsilon_0$ sufficiently small, depending on $\|\Phi\|_{\mathcal{L}_2^{0,1}}$, T (so that $\epsilon_0^2T$ is small) and $\alpha$, and for $\epsilon<\epsilon_0$,
\begin{align}
\label{st1}
\mathbb{P}\left(|c^{\epsilon}\left(\tau\right)-c_0|^2\ge\alpha^2\right)&
\le\mathbb{P}\left(\epsilon\left|\int_0^{\tau}\left(u^{\epsilon,\varphi_{c_0}}(s),dW(s)\right)_{\xLtwo}\right|
\ge\frac{\alpha}{4}\right)\nonumber\\
&\le\mathbb{P}\left(\epsilon\sup_{t\in[0,T]}\left|\int_0^{t}
\left(u^{\epsilon,\varphi_{c_0},\tau^{\epsilon}_{\alpha}}(s),dW(s)\right)_{\xLtwo}\right|
\ge\frac{\alpha}{4}\right).
\end{align}
where $u^{\epsilon,\varphi_{c_0},\tau^{\epsilon}_{\alpha}}$ is the process stopped 
at time $\tau^{\epsilon}_{\alpha}$.
Similarly, using as well the following property
\begin{itemize}
\item[$(P)\;$] When $|c^{\epsilon}\left(\tau\right)-c_0|\le\alpha\le\alpha_0$, then for
some $C$ depending only on $c_0$ and $\alpha_0$,
$\|\varphi_{c^{\epsilon}(t)}-\varphi_{c_0}\|_{\xHone}\le C\alpha$, for all $t\le \tau$,
\end{itemize}
we obtain that there exists $\epsilon_0$ sufficiently small, depending on $\|\Phi\|_{\mathcal{L}_2^{0,1}}$, T (with $\epsilon_0^2T$ small) and $\alpha$, such that
for all $\epsilon<\epsilon_0$,
\begin{align}
\label{st2}
&\mathbb{P}\left(\|\epsilon\eta^{\epsilon}\left(\tau\right)\|_{\xHone}^2\ge\alpha^2\right)\nonumber\\
\le &  \mathbb{P}\left(\epsilon\sup_{t\in[0,T]}\left|\int_0^{t}
\left(u^{\epsilon,\varphi_{c_0},\tau_{\alpha}^{\epsilon}}(s),dW(s)\right)_{\xLtwo}\right|
\ge\frac12\sqrt{\frac{\alpha^2}{10}}\right)\nonumber\\
+ &\mathbb{P}\left(\epsilon\sup_{t\in[0,T]}\left|\int_0^{t}
\left(\partial_xu^{\epsilon,\varphi_{c_0}\tau_{\alpha}^{\epsilon}}(s),\partial_xdW(s)\right)_{\xLtwo}\right|
\ge\frac{\alpha^2}{10}\right)\nonumber\\
+ &\mathbb{P}\left(\epsilon\sup_{t\in[0,T]}\left|\int_0^{t}
\left(\left(u^{\epsilon,\varphi_{c_0}\tau_{\alpha}^{\epsilon}}(s)\right)^2,dW(s)\right)_{\xLtwo}\right|
\ge\frac{\alpha^2}{10}\right)\nonumber\\
+ &\mathbb{P}\left(\epsilon\sup_{t\in[0,T]}\left|\int_0^{t}
\left(u^{\epsilon,\varphi_{c_0},\tau_{\alpha}^{\epsilon}}(s),dW(s)\right)_{\xLtwo}\right|
\ge\frac{\alpha^2}{10c_0}\right).
\end{align}
Let us denote by $Z_i(t)$ for $i=\{1,2,3\}$ the stochastic
integrals arising in the right hand sides of (\ref{st1}) and (\ref{st2}).
We obtain exponential tail estimates for each of the
above probabilities in the usual way, see for example \cite{P2}, Theorem
2.1. We introduce the function $f_l(x)=\sqrt{1+lx^2}$, where $l$ is
a positive parameter. We then apply the It\^o formula to
$f_l\left(Z_i(t)\right)$ and each process decomposes into
$1+E_{l,i}(t)+R_{l,i}(t)$ where
\begin{equation*}
E_{l,i}(t)=\int_0^t\frac{lZ_i(t)}{\sqrt{1+lZ_i(t)^2}}dZ_i(t)
-\frac12\int_0^t\left(\frac{lZ_i(t)}{\sqrt{1+lZ_i(t)^2}}\right)^2d<Z_i>_t,
\end{equation*}
and
\begin{equation*}
R_{l,i}(t)=\frac12\int_0^t\left(\frac{lZ_i(t)}{\sqrt{1+lZ_i(t)^2}}\right)^2d<Z_i>_t
+\frac12 \int_0^t\frac{d<Z_i>_t}{\left(1+lZ_i(t)^2\right)^{3/2}}.
\end{equation*}
Let us for example consider $Z_2$. Given
$\left(e_j\right)_{j\in\xN}$ a complete orthonormal system of
$\xLtwo$,
\begin{equation*}
<Z_2>_t=\int_0^t
\sum_{j\in\xN}\left(\partial_xu^{\epsilon,\varphi_{c_0}\tau_{\alpha}^{\epsilon}},\partial_x\Phi
e_j\right)_{\xLtwo}^2(s)ds; 
\end{equation*} 
thus, using the H\"older inequality and the property $(P)$ we have,
for some constant $C(\alpha,c_0)=C\alpha^2+\|\partial_x\varphi_{c_0}\|_{\xHone}$, and any $t$,
$$<Z_2>_t\le C(\alpha,c_0)\|\Phi^*\|_{\mathcal{L}^{-1,0}}^2t.$$
Then $$|R_{l,2}(t)|\le lC(\alpha,c_0)\|\Phi\|_{\mathcal{L}^{0,1}}^2T.$$
The same bound also holds for $Z_1$ and $Z_3$.
We may thus write for any $i$ and constants
$\delta_i>0$,
\begin{align*}
& \ \mathbb{P}\left(\sup_{t\in[0,T]}\left|Z_i(t)\right|\geq
\frac{\delta_i}{\epsilon}\right)\\
= & \ \mathbb{P}
\left(\sup_{t\in[0,T]}\exp\left(f_l(Z_i(t))
\right)\geq\exp\left(f_l\left(\frac{\delta_i}{\epsilon}\right)\right)\right)\\
\leq & \ \mathbb{P}\left(\sup_{t\in[0,T]}\exp\left(E_{l,i}(t)
\right)\geq\exp\left(f_l\left(\frac{\delta_i}{\epsilon}\right)
-1-l C(\alpha,c_0)\|\Phi\|_{\mathcal{L}^{0,1}}^2T\right)\right).
\end{align*}
The Novikov condition is also satisfied and $E_{l,i}(t)$ is such
that $\left(\exp\left(E_{l,i}(t)\right)\right)_{t\ge0}$ is a
uniformly integrable martingale. The Doob inequality then gives
\begin{align*}
& \ \mathbb{P}\left(\sup_{t\in[0,T]}\exp\left(E_{l,i}(t)
\right)\ge\exp\left(f_l\left(\frac{\delta_i}{\epsilon}\right)
-1-lC(\alpha,c_0)\|\Phi\|_{\mathcal{L}^{0,1}}^2T\right)\right)\\
\le & \ \exp\left(-f_l\left(\frac{\delta_i}{\epsilon}\right)
+1+lC(\alpha,c_0)\|\Phi\|_{\mathcal{L}^{0,1}}^2T\right)\mathbb{E}\left[\exp\left(E_{l,i}(T)\right)\right].
\end{align*}
Since
$\exp\left(E_{l,i}(T)\right)$ is an exponential martingale $\mathbb{E}\left[\exp\left(E_{l,i}(T)\right)\right]=\mathbb{E}\left[\exp\left(E_{l,i}(0)\right)\right]=1$.
For $\epsilon_0$ small enough we have for $\epsilon<\epsilon_0$
$$C(\alpha,c_0)\|\Phi\|_{\mathcal{L}^{0,1}}^2T<\frac{\delta_i^2}{2\epsilon^2}$$
which implies that the $l$-derivative at 0 of the function in the exponential bound is negative.
Then, optimizing on the parameter $l$, we obtain the minimum in $l$ of the upper bound,
which has the form
$$\exp\left(1-\frac{C(\alpha,c_0)}{\epsilon^2 T\|\Phi\|_{\mathcal{L}^{0,1}}^2}\right),$$
with possibly another constant $C(\alpha,c_0)$.
Using the largest of all constants in the exponentials for each
tail probabilities, and taking a constant slightly bigger and $\epsilon_0$ smaller if
necessary, the multiplicative constant in front of the exponential decay disappears
and the result follows.\hfill $\square$\\

\begin{rmrk}
As will appear elsewhere, Theorem \ref{tdBD} also holds for a noise of multiplicative
type, and an exponential upper bound holds as well.
\end{rmrk}

If we now consider the sequence of operators $\Phi_n$ mimicking the spatially homogeneous noise with
covariance operator $(I-\Delta)^{-1}$, denoting
the exit times $\tau_{\alpha}^{n,\epsilon}$, we obtain
the following statement:

\begin{crllr}\label{c1}
For $T>0$, $0<\alpha<\alpha_0$ and $n$ fixed,
there exists a constant $C(\alpha,c_0)$, depending on $\alpha$ and $c_0$, and
there exists $\epsilon_0>0$ with $\epsilon_0^2T$ sufficiently small with respect to
$\|\Phi_n\|_{\mathcal{L}_2^{0,1}}$ and $\alpha$, such that for every
$\epsilon<\epsilon_0$,
\begin{equation}\label{e6}
\mathbb{P}\left(\tau_{\alpha}^{n,\epsilon}\le
T\right)\le\exp\left(-\frac{C(\alpha,c_0)}{\epsilon^2T}\right).
\end{equation}
In particular, we have the following double asymptotic result
$$\overline{\lim}_{n\rightarrow\infty}\overline{\lim}_{\epsilon\rightarrow0}\epsilon^2\log\mathbb{P}\left(\tau_{\alpha}^{n,\epsilon}\le
T\right)\le-\frac{C(\alpha,c_0)}{T}.$$
\end{crllr}

Let us now prove that when considering the sequence of operators $\Phi_n$ as above,
one can obtain a lower bound of the
same order in the parameters $T$ and $\epsilon$ as the upper bound in Corollary \ref{c1}. We make use of the
approximation via a modulated soliton together with the LDP obtained in Theorem \ref{t1},
and again minimize the rate function on a smaller set of controls giving rise to a set of parameterized 
exit paths.

\begin{prpstn}\label{p2}
For every $T$ and $\alpha$ positive, there exists a
constant $C(\alpha, c_0)$ which depends on $c_0$ and $\alpha$ but not on $T$, such that
$$\underline{\lim}_{n\to\infty}\underline{\lim}_{\epsilon\to0}\epsilon^2\log\mathbb{P}
\left(\tau_{\alpha}^{n,\epsilon}\le T\right)\ge-\frac{C(\alpha,c_0)}{T}.$$
\end{prpstn}

\noindent
{\bf Proof.}
Let us denote by $\mathcal{U}_{\alpha_0}$ the open set
$$\mathcal{U}_{\alpha_0}=\left\{\varphi_c(\cdot-y)+g,\ g\in\xHone:\
\|g\|_{\xHone}<\alpha_0,\ y\in\xR,\ |c-c_0|< \alpha_0\right\}.$$
We know from \cite{dBD2} that the velocity is obtained via a
continuous mapping
$\mathcal{C}$ from $\mathcal{U}_{\alpha_0}$ to $\xR$
such that
$$c^{\epsilon}(t)=\mathcal{C}\left(u^{n,\epsilon,\varphi_{c_0}}(t)\right).$$
Also for fixed $n$, we write that for $0<2\alpha<\alpha_0$,
$$\mathbb{P}\left(\tau_{\alpha}^{n,\epsilon}\le T\right)\ge\mathbb{P}
\left(u^{n,\epsilon,\varphi_{c_0}}\in B\right)$$ where
$$B=\left\{u\in\xCn([0,T];\xHone):\ \forall t\in[0,T],\ u(t)\in\mathcal{U}_{\alpha_0},\
|\mathcal{C}\left(u(T)\right)-c_0|>\alpha\right\}.$$
Theorem \ref{t1} then leads to the following lower bound for
$\underline{\lim}_{\epsilon\to0}\epsilon\log\mathbb{P}\left(\tau_{\alpha}^{n,\epsilon}\le
T\right)$:
\begin{equation}
\label{inf1}
-\frac12\inf\left\{\|h\|_{\xLtwo(0,T;\xLtwo)}^2,\ h:\ \forall t\in[0,T],\ \mathbf{S}^{n,\varphi_{c_0}}(h)(t)\in
\mathcal{U}_{\alpha_0},\ \left|\mathcal{C}
\left(\mathbf{S}^{n,\varphi_{c_0}}(h)(T)\right)-c_0\right|>\alpha\right\}.
\end{equation}
Let us, as in the proof of Proposition \ref{p3}, replace in a first step the above variational problem
by a variational problem for $\Phi=\Phi_{H}=(I-\Delta)^{-1/2}$ and $\alpha$ by $2\alpha$.
Minimizing on a smaller set, we obtain
$$
-\frac12\inf\left\{\|h\|_{\xLtwo(0,T;\xLtwo)}^2,\ h:\ \forall t\in[0,T],\ \mathbf{S}^{H,\varphi_{c_0}}(h)(t)\in
\mathcal{U}_{\alpha_0},\ \mathcal{C}\left(\mathbf{S}^{H,\varphi_{c_0}}(h)(T)\right)-c_0=3\alpha
\right\}.$$
We minimize on an even smaller set, considering solutions of the
controlled equation which are modulated solitons of the form
$$\varphi_{c(t)}\left(x-\int_0^{t}c(s)\right)$$
where the one dimensional paths $c$ are assumed to belong to
$\xCn^1([0,T]; (0,+\infty))$. The boundary conditions are thus that $c(0)=c_0$ and
$c(T)=c_0+2\alpha$. A control $h_c$ associated to such
a solution is given by
\begin{align*}
h_c(t,x)&=\Phi_{H}^{-1}\left(\partial_t\mathbf{S}^{H,\varphi_{c_0}}(h_c)+\partial_x^3\mathbf{S}^{H,\varphi_{c_0}}(h_c)
+\partial_x\left(\mathbf{S}^{H,\varphi_{c_0}}(h_c)^2\right)\right)(t,x)\\
&=c'(t)\left.\Phi_{H}^{-1}\partial_c\varphi_c\right|_{c=c(t)}\left(x-\int_0^tc(s)ds\right).
\end{align*}
Note that this control is the same as in Proposition \ref{p3}, only the terminal
boundary condition changes.
We thus obtain the lower bound
\begin{align}
\label{inf2}
-\frac12\inf&\left\{\int_0^T\left\|c'(t)(I-\Delta)^{1/2}
\left.\partial_c\varphi_c\right|_{c=c(t)}
\left(x-\int_0^tc(s)ds\right)\right\|_{\xLtwo}^2dt,\right.\nonumber\\
&\quad\left.c\in\xCn^1([0,T]; (0,+\infty)):\ c(0)=c_0,\ c(T)=c_0+2\alpha\right\}.
\end{align}
We now have to solve a problem of the calculus of variations. Our aim
is to find the optimal paths $c$ among the set of constrained paths
minimizing the path integral. The integral may be written, with the same
function $g$ as in Proposition \ref{p3}, as
$$\int_0^T\left(c'(t)\right)^2g\left(c(t)\right)dt.$$
Using successively the change of variables $t=Tu$ and the change of
unknown function $v(u)=c(Tu)$, we obtain an upper bound of the form
$$\frac1T\int_0^1\left(v'(u)\right)^2g(v(u))du$$
for functions $v$ which are $\xCn^1([0,T];(0,+\infty))$ and satisfy the two
boundary conditions $v(0)=c_0$ and $v(1)=c_0+2\alpha$ independent of $T$.
We recall that \linebreak $g(c)=\|(I-\Delta)^{1/2} \partial_c \varphi_c\|_{\xLtwo}^2$ so that
$s \mapsto g(v(s))$ is bounded on $[0,1]$ for any $v\in \xCn^1([0,T];(0,+\infty))$.
Hence we deduce that the infimum in (\ref{inf2}) is bounded above by $\frac{C(\alpha,c_0)}{T}$.
Now, coming back to the case where $\Phi_n$ is the square root of the covariance operator of the noise,
we start from a path $c$ obtained from $v$ which say minimizes the objective function
in the above problem of the calculus of variations, though following the above argument it does
not really matter.
Then, the control $h_c$ is such that
$$\mathcal{C}\left(\mathbf{S}^{H,\varphi_{c_0}}(h_c)(T)\right)=c_0+2\alpha.$$
From the continuity of the mild-solution of the control map
with respect to the convolution of the semi-group
with the control, used to prove the LDP, and the continuity
of this last quantity with respect to the control, and using also the continuity of
$\mathcal{C}$ with respect to $u \in \mathcal{U}_{\alpha_0}$, we know that for sufficiently large $n$
\begin{align*}
&\mathcal{C}\left(\mathbf{S}^{n,\varphi_{c_0}}(h_c)(T)\right)\\
= \ &\mathcal{C}\left(\mathbf{S}^{H,\varphi_{c_0}}\left((I-\Delta)^{1/2}
(I-\Delta+\frac{1}{n}(x^2I-\Delta)^k)^{-1/2}h_c\right)(T)\right)\\
> \ &\alpha.
\end{align*}
We deduce that the inf-limit as $n$ goes to infinity of the infimum in (\ref{inf1}) is again
bounded above by $\frac{C(\alpha, c_0)}{T}$ and 
this ends the proof of Proposition \ref{p2}. \hfill $\square$\\

As a consequence of our two bounds, the typical time scale on which the solution remains in the
neighborhood of the modulated soliton is indeed $1/\epsilon^2$.\vspace{0.5cm}

\noindent {\bf Acknowledgments.} The authors are grateful to the
Centro di Ricerca Matematica Ennio De Giorgi in Pisa for allowing us
to participate to the spring 2006 research period when this research was initiated.

\end{document}